\documentclass[12pt]{amsart}
\usepackage{amsmath}
\usepackage{amsthm}
\usepackage{amssymb}
\usepackage{amsfonts,mathrsfs}
\usepackage{amsxtra}
\usepackage{epsfig}
\usepackage{verbatim}
\usepackage{color}
\usepackage{hyperref}

\theoremstyle{plain}

\newtheorem{lemma}[equation]{Lemma}

\theoremstyle{remark}

\include{header}

\begin{document}

\title[Proofs in Number Theory and Algebra Through Bergman Spaces]{The Fundamental Theorem of Algebra and the Divergence of Reciprocals of Primes looked at through Bergman Spaces}%{A Look at Elementary Number Theory Through Bergman Spaces}
\author{Yunus E. Zeytuncu}
\email{zeytuncu@math.tamu.edu}
\begin{abstract}
Two well known facts from elementary number theory are proven by using Bergman spaces.
\end{abstract}
%\keywords{32A36, Bergman Spaces, Fundamental Theorem of Algebra, Harmonic series of primes}
\maketitle
\section{Introduction}

In this note, we present proofs for two well known facts in elementary number
theory through Bergman spaces terminology. Both of these facts have complex analytic proofs (in addition to proofs by other techniques); therefore, it is not
surprising to be able to write proofs by using Bergman spaces. However, we present these proofs here with hope for more connections between Bergman spaces and number theory. 

In the second section, before we talk about the proofs, we start with a brief introduction to Bergman spaces.

In the third section, we present a proof of The Fundamental Theorem of Algebra (FTA). We refer to \cite{Fine97} for the history and six different proofs 
(including a complex analysis proof) of the FTA. 
There are more than 20 articles in MAA's Monthly that present other proofs. The most recent one, \cite{Lazer}, claims that at least 80 proofs exist in 
the literature and the note, \cite{What}, can be seen for comments on some of the proofs. 

In the fourth section, we look at the harmonic series of prime numbers. It is a well-known fact that the series of reciprocals of prime numbers diverges. 
This was first noted by Euler in 18th century. 
Section 33 of \cite{Euler1748} and the references therein can be consulted to see \textit{how Euler did it}. 
Two other rigorous proofs can be additionally found in the first section of \cite{Ziegler99}  and on page 297 of \cite{Nowak2000}. 
The one in \cite{Ziegler99} is due to Erd\"os. In the fourth section, we present another proof of this fact.\\

\section{Bergman Spaces over Discs}

Bergman spaces are defined for any subset of the complex plane $\mathbb{C}$ but for simplicity we will only discuss Bergman spaces of discs.
For any $R>0$, let $D_R$ denote the open disc $|z|<R$ in $\mathbb{C}$ (when $R=1$, $\mathbb{D}$ is used to denote the unit disc) and
$L^2(D_R)$ denote the space of square integrable functions on $D_R$. $L^2(D_R)$ is a Hilbert space with the natural inner product structure given by
$$\left<f,g\right>=\int_{D_R}f(z)\overline{g(z)}dA(z)$$ for $f$ and $g$ in $L^2(D_R)$. Here, $dA(z)$ stands for the Lebesgue measure on $\mathbb{C}$.
The subspace of holomorphic functions that are in $L^2(D_R)$ is called the Bergman space of $D_R$ and denoted by $A^2(D_R)$. 
It is an important consequence of the Cauchy integral formula and the Cauchy-Schwarz inequality that $A^2(D_R)$ is a closed subspace of $L^2(D_R)$ 
(see \cite{Duren04} or \cite{KrantzGFT06}). Therefore, there is the following orthogonal decomposition:

\begin{equation}\label{decomposition}
L^2(D_R)=A^2(D_R)\oplus A^2(D_R)^{\perp}.
\end{equation}

The Taylor series expansions of holomorphic functions give a direct way of calculating inner products and norms on $A^2(D_R)$. 
Namely, for two functions $f$ and $g$ in $A^2(D_R)$ with $f(z)=\sum_{n=0}^{\infty}f_nz^n$ and $g(z)=\sum_{n=0}^{\infty}g_nz^n$ a routine calculation gives
\begin{align}\label{identity}
\left<f,g\right>=\pi\sum_{n=0}^{\infty}R^{2n+2}\frac{f_n\overline{g_n}}{n+1}~\text{ and }~
||f||^2&=\pi \sum_{n=0}^{\infty}R^{2n+2}\frac{|f_n|^2}{n+1}~.
\end{align}
\vskip .5cm
\section{The Fundamental Theorem of Algebra}

\allowdisplaybreaks{
Let $P(z)=a_nz^n+\dots+a_1z+a_0$ be a non-constant complex polynomial, where $a_n \not = 0, a_0\not=0$ and $n\geq 2$ (without these restriction the FTA follows trivially). 
Suppose that $P(z)$ has no zeros in $\mathbb{C}$. Then $F(z)=\overline{\frac{1}{P(z)}}$ is a well-defined continuous function on $\mathbb{C}$. 
The triangle inequality implies that there exists $R_0>0$ such that, for any $|z|>R_0$,
\begin{align*}
|F(z)|=\frac{1}{|P(z)|}&=\frac{1}{|a_nz^n|\left|1+\frac{a_{n-1}}{a_n}\frac{1}{z}+\dots+\frac{a_0}{a_n}\frac{1}{z^n}\right|}\\
&\leq\frac{1}{|a_nz^n|\left(1-\left|\frac{a_{n-1}}{a_n}\frac{1}{z}+\dots+\frac{a_0}{a_n}\frac{1}{z^n}\right|\right)}\\
&\leq \frac{2}{|a_nz^n|}.
\end{align*}
This, in particular, implies that for any $R>R_0$,
\begin{align*}
\int_{R_0<|z|<R}|F(z)|^2dA(z)&\leq \int_{R_0 < |z|<R}\frac{4}{|a_n|^2|z|^{2n}}dA(z)= \frac{8\pi}{|a_n|^2}\int_{R_0}^R\frac{1}{r^{2n}}rdr\\
&= \frac{4\pi}{|a_n|^2(n-1)}\left(\frac{1}{R_0^{2n-2}}-\frac{1}{R^{2n-2}}\right)\\
&\leq \frac{4\pi}{R_0^{2n-2}|a_n|^2(n-1)}.
\end{align*}}

\noindent As a consequence of this estimate, there exists a constant $M>0$ such that for any $R>0$,
\begin{equation}\label{est}
\int_{|z|<R}|F(z)|^2dA(z)\leq M.
\end{equation}

By the estimate \eqref{est}, we get $F\in L^2(D_R)$ for any $R>0$. The decomposition \eqref{decomposition} gives that for any $R>0$ there exist functions $G_R$ and $H_R$ such that $G_R \in A^2(D_R)$, $H_R\in A^2(D_R)^{\perp}$ and $F=G_R+H_R$. The key observation is the following lemma.
\begin{lemma}\label{lemma}
For any $R>0$, 
\begin{equation*}
G_R(z)=\overline{\frac{1}{a_0}}~.
\end{equation*}
\end{lemma}

\noindent In more technical terminology, the lemma says that on any domain $D_R$ the \textit{Bergman projection} of $F(z)$ is the constant $\overline{\frac{1}{a_0}}$. 

\begin{proof}
By the assumption, the function $\frac{1}{P(z)}$ is entire. Therefore, it has the Taylor series expansion, $\frac{1}{P(z)}=\sum_{n=0}^{\infty}b_nz^n$, that converges uniformly on any $D_R$ and it is clear that $b_0=\frac{1}{P(0)}=\frac{1}{a_0}$.\\

\noindent Using this representation, we get
\begin{equation*}
F(z)=\overline{\frac{1}{a_0}}+\sum_{n=1}^{\infty}\overline{b_nz^n}.
\end{equation*}

\noindent It is simple to check that for any $n\geq1$, the monomial $\overline{z}^n\in A^2(D_R)^{\perp}$ for any $R>0$. This concludes the proof of the lemma.
\end{proof}

Now, for any $R>0$, the Pythagorean theorem and Lemma \ref{lemma} together imply that
\begin{align}\label{lowest}
\int_{D_R}|F(z)|^2dA(z)\geq \int_{D_R}|G_R(z)|^2dA(z)
=\int_{D_R}\frac{1}{|a_0|^2}dA(z)
=\frac{\pi R^2}{|a_0|^2}~.
\end{align}

Finally, the combination of the inequalities \eqref{est} and \eqref{lowest} gives
\begin{equation*}
M\geq \int_{D_R}|F(z)|^2dA(z) \geq \frac{\pi R^2}{|a_0|^2}
\end{equation*}
for any $R>0$. This is a contradiction so $P(z)$ has a zero in $\mathbb{C}$.
\vskip .5cm

\section{The Harmonic Series of Prime Numbers}

Let $P$ denote the set of prime numbers i.e. $P=\{2,3,5,7,\dots\}$. We look at the function $\mathcal{P}(z)$ defined by $$\mathcal{P}(z)=\sum_{p\in P}z^p.$$ 
It is clear that $\mathcal{P}(z)$ is holomorphic on $\mathbb{D}$ and the identity \eqref{identity} indicates
\begin{equation*}
||\mathcal{P}||^2=\pi \sum_{p\in P}~\frac{1}{p+1}\approx \sum_{p\in P}~\frac{1}{p}~.
\end{equation*}

Therefore $\mathcal{P}(z)$ is square integrable on $\mathbb{D}$ if and only if $\sum_{p\in P}~\frac{1}{p}$ is convergent. 
So in order to prove that the harmonic series of prime numbers diverges;  it is enough to prove that $\mathcal{P}(z) \not\in A^2(\mathbb{D})$.

The proof is by \textit{reductio ad absurdum}. Namely, we assume that $\mathcal{P}(z)\in A^2(\mathbb{D})$ and we show that this implies $\frac{1}{1-z}\in A^2(\mathbb{D})$. 
However, this is a contradiction since $\frac{1}{1-z}$ is \textit{too} singular at $z=1$ to be square integrable on $\mathbb{D}$. 
A direct computation shows that for $t\in\mathbb{R}$; the holomorphic function $(1-z)^t$ is square integrable on $\mathbb{D}$ if and only if $t>-1$.

First note that if $g(z)=f(z^m)$ for some $f\in A^2(\mathbb{D})$ such that $f(0)=0$ and $m\in \mathbb{N}$ then $g$ is also in $A^2(\mathbb{D})$ and \eqref{identity} implies

\begin{equation}\label{comparison}
||g||^2<~\frac{2}{m}||f||^2.
\end{equation}

As mentioned above, now let us suppose that $\mathcal{P}(z)\in A^2(\mathbb{D})$. This implies that
$\sum_{p\in P}\frac{1}{p}$ is convergent and there exists a prime number $p_k$ such that
\begin{equation}\label{tail}
\sum_{p\in P,~ p\geq p_k}\frac{1}{p}<1.
\end{equation}
Let $P_1$ be the set of prime numbers less than $p_k$
and $P_2$ be the remaining ones i.e. $P_1=\{2,3,\dots,p_{k-1} \}$ and $P_2=\{p_k,p_{k+1},\dots \}$.
Also consider two subsets of the natural numbers $N_1$ and $N_2$ that are defined by:
\begin{align*}
N_1&=\{n\in \mathbb{N}~|~\text{ prime decomposition of }n \text{ contains primes only from }P_1\}\\
N_2&=\{n\in \mathbb{N}~|~\text{ prime decomposition of }n \text{ contains primes only from }P_2\}.\\
\end{align*}

\textit{Step One.} We define another function $F(z)=\sum_{n\in N_2}z^n$. At the end of this section we verify that the assumption above guarantees $F$ is in $A^2(\mathbb{D}).$   
For the moment, we take this granted and we use $F(z)$ to decompose $\frac{1}{1-z}$ as sum of mutually orthogonal functions. Namely, we show that

\begin{align}
\frac{1}{1-z}&=\sum_{n=0}^{\infty}z^n
=1+z+F(z)+\sum_{k\in N_1}\left(z^k+F(z^k)\right).
\end{align}

It is clear that $F(z^k)$ is a sum of monomials for any $k\in N_1$ and if we also consider $F(z^k)$ as a list of monomials then any monomial $z^n$ appears somewhere in the set $\{z^k,F(z^k)\}_{k\in N_1}\bigcup\{1,z,F(z)\}$ since any integer is product of an integer from $N_1$ and an integer from $N_2$. Moreover any monomial appears in this list only once. Therefore any pair of functions in $\{z^k,F(z^k)\}_{k\in N_1}\bigcup\{1,z,F(z)\}$ have no monomial in common so they are mutually orthogonal. By using the orthogonality and the estimate \eqref{comparison} we obtain:

\begin{align*}
||\frac{1}{1-z}||^2&=||1||^2+||z||^2+||F(z)||^2+\sum_{k\in N_1}\left(||z^k||^2+||F(z^k)||^2\right)\\
&\leq \frac{3\pi}{2}+||F(z)||^2+\sum_{k\in N_1}\left(\frac{\pi}{k+1}+\frac{2}{k}||F(z)||^2\right)\text{ here we use the estimate }\eqref{comparison}\\
&\leq \frac{3\pi}{2}+||F(z)||^2+\left(\pi+2||F(z)||^2\right)\sum_{k\in N_1}\frac{1}{k}
\end{align*}

When we rewrite the last sum as:
\begin{align*}
\sum_{k\in N_1}\frac{1}{k}&=\left(1+\frac{1}{2}+\frac{1}{2^2}+\dots\right)\left(1+\frac{1}{3}+\frac{1}{3^2}+\dots\right)
\dots\left(1+\frac{1}{p_{k-1}}+\frac{1}{p_{k-1}^2}+\dots\right)\\
&=\left(\frac{1}{1-\frac{1}{2}}\right)\left(\frac{1}{1-\frac{1}{3}}\right)\dots\left(\frac{1}{1-\frac{1}{p_{k-1}}}\right)
\end{align*}
we note that it is finite and conclude that $\frac{1}{1-z}\in A^2(\mathbb{D})$. This is the desired contradiction to finish the proof.\\

\textit{Step Two.} It now remains to verify that the assumption \eqref{tail} indeed implies that $F(z)\in A^2(\mathbb{D})$. For this purpose, we first look at the function $\mathcal{Q}(z)=\sum_{p\in P_2}z^p$. It is evident that $\mathcal{Q}\in A^2(\mathbb{D})$ . We consider the functions $H_l(z)=\mathcal{Q}(z^l)$ for any $l\in N_2$. We put the numbers in $N_2$ in the increasing order and label them as $l_1,l_2,\dots$. After this we modify each function $H_l$ to a new function $G_l$ as follows:
\begin{align*}
G_{l_1}(z)&=H_{l_1}(z)\\
G_{l_2}(z)&=H_{l_2}(z) \text{ but drop the monomials which appear in } G_{l_1}\\
&~\vdots\\
G_{l_i}(z)&=H_{l_i}(z) \text{ but drop the monomials which appear in } G_{l_1},\dots,G_{l_{i-1}}\\
&~\vdots
\end{align*}

This modification guarantees that $G_l$'s and $\mathcal{Q}$ do not have any monomial in common and  all of $G_l$'s and $\mathcal{Q}$ are orthogonal to each other. 
Furthermore, any number in $N_2$ appears as a power either in $\mathcal{Q}(z)$ or in one of $G_l(z)$'s therefore we have the following orthogonal decomposition of $F$

\begin{equation}
F(z)=\mathcal{Q}(z)+\sum_{l\in N_2}G_l(z).
\end{equation}

The identity \eqref{identity} says $||G_l||\leq||H_l||$ for any $l \in N_2$ and the inequality \eqref{comparison} (note that $\mathcal{Q}(0)=0$) gives 

\begin{equation}\label{estimate}
||G_l||^2\leq||H_l||^2\leq\frac{2}{l}||\mathcal{Q}||^2 ~\text{ for any } l \in N_2.
\end{equation}

Hence we get
\begin{align*}
||F||^2&=||\mathcal{Q}||^2+\sum_{l\in N_2}||G_l||^2~\text{ by orthogonality}\\
&\leq ||\mathcal{Q}||^2+\sum_{l\in N_2}\frac{2}{l}||\mathcal{Q}||^2~\text{ by }\eqref{estimate}\\
&\leq 2||\mathcal{Q}||^2\left(1+\sum_{l\in N_2}\frac{1}{l}\right)
\end{align*}

We can rewrite the last series as follows:
\begin{align*}
\sum_{l\in N_2}\frac{1}{l}&=\sum_{p\in P_2}\frac{1}{p}+\sum_{p,q\in P_2}\frac{1}{pq}+\sum_{p,q,r\in P_2}\frac{1}{pqr}+\dots\\
&\leq\left(\sum_{p\in P_2}\frac{1}{p}\right)+\left(\sum_{p\in P_2}\frac{1}{p}\right)^2+\left(\sum_{p\in P_2}\frac{1}{p}\right)^3+\dots
\end{align*}
Since $\sum_{p\in P_2}\frac{1}{p}<1$ we get that $\sum_{l\in N_2}\frac{1}{l}$ is convergent and we indeed verify that $F$ is in $A^2(\mathbb{D})$.

\section{Remarks}

Some other well known facts from elementary number theory also can be rephrased in the terminology of Bergman spaces. Two examples of this translation are the followings.

\noindent \textbf{Example 1.} Let $g_N(z)=\sum_{k=N+1}^{2N}z^k$ and $q_N(z)=\sum_{p\in P, ~p\leq 2N}z^p$ be two polynomials. The identity \eqref{identity} implies that
\begin{align*}
\left<g_N,q_N\right>=\pi\sum_{p\in P\cap\{N+1,\dots,2N\}}\frac{1}{p+1}.
\end{align*}
Therefore $\left<g_N,q_N\right>$ is non-zero if and only if there is at least one prime in the set $\{N+1,\dots,2N\}$. The latter is exactly \textit{Bertrand's postulate} 
(see \cite{Ziegler99}) so if one can directly prove that $\left<g_N,q_N\right>$ is non-zero then we get an alternative proof of Bertrand's postulate.

\noindent \textbf{ Example 2.} Let $P_t$ denote the set of twin primes i.e. $P_t=\{p\in P:~p+2\in P\}$ and let
$\mathcal{T}(z)=\sum_{p\in P_t}z^p$. The identity \eqref{identity} implies that
\begin{align*}
||\mathcal{T}||^2=\pi \sum_{p\in P_t}~\frac{1}{p+1}\approx \sum_{p\in P_t}~\frac{1}{p}~.
\end{align*}
Therefore $\mathcal{T}$ is in $A^2(\mathbb{D})$ if and only if the sum of reciprocals of twin primes converges. 
The latter is exactly \textit{Brun's Theorem} (see \cite{Leveque77}) so if one can directly prove that $\mathcal{T}$ is in $A^2(\mathbb{D})$ 
then we get an alternative proof of Brun's Theorem.

\section*{Acknowledgments} I'd like to thank T. Kaptanoglu and D. Vukotic for useful comments on this paper. I also thank the anonymous referee for helpful recommendations to improve the presentation of the paper.

\bibliographystyle{plain}
\bibliography{PrimesBergman}
\vskip 1cm
\end{document}